\documentclass[a4paper,12pt]{article}

\usepackage{amsthm,amsmath,amssymb}
\usepackage{graphpap,latexsym}

\usepackage{epsfig}

\usepackage{color}

\definecolor{blue3}{rgb}{0,.3,1}
\definecolor{green1}{rgb}{0.1,0.65,0}
\definecolor{brown}{rgb}{.4,.25,0.05}

\numberwithin{equation}{section}

\newtheorem{theorem}{Theorem}[section]

\newtheorem{remark}[theorem]{Remark}

\newtheorem{hypo}[theorem]{Hypothesis}

\newcommand{\behy}{\begin{hypo}\begin{sl}}
\newcommand{\enhy}{\end{sl}\end{hypo}}

\def\bfsi{{\boldsymbol \sigma}}
\def\bfe{{\boldsymbol \varepsilon}}

\def\bfde{{\boldsymbol \delta}}
\def\vr{\varrho}

\newcommand{\real}{\mathbb{R}}

\newcommand{\nat}{\mathbb{N}}

\def\vp{\varphi}
\def\e{\varepsilon}

\def\dive{\mbox{\rm div\,}}

\newcommand{\io}{\int_\Omega}
\newcommand{\ipo}{\int_{\partial\Omega}}

\def\itt{\int_0^t}

\def\bfb{\mathbf{b}}

\def\bff{\mathbf{f}}
\def\bfg{\mathbf{g}}
\def\bfk{\mathbf{k}}
\def\bfn{\mathbf{n}}
\def\bfq{\mathbf{q}}
\def\bfu{\mathbf{u}}
\def\bfv{\mathbf{v}}

\def\bfdd{\mbox{\bf \,:\,}}

\def\no{\nonumber}
\def\dd{\mathrm{\,d}}

\def\AA{\mathcal{A}}
\def\BB{\mathcal{B}}
\def\CC{\mathcal{C}}

\def\tetak{\theta_k}
\def\tetakm{\theta_{k-1}}
\def\Uk{U_k}
\def\Ukm{U_{k-1}}
\def\chik{\chi_k}
\def\chikm{\chi_{k-1}}

\def\bdot#1{\hspace{-1pt}\stackrel{\bullet}{#1}\!\!}
\def\dUk{\bdot{U}_k}
\def\dchik{\bdot{\chi}_k}
\def\Uko{U_{k\Omega}}

\def\teta{\theta}
\makeindex             % used for the subject index
\begin{document}

\title{Well-posedness of an extended model\\ for water-ice phase transitions\thanks{Supported by GA\v CR
Grant No. P201/10/2315. A large part of this work was done during
the visits of ER in Prague and of PK in Milan and Pavia.}}

\author{Pavel Krej\v{c}\'{i}
\thanks{Institute of Mathematics,
Czech Academy of Sciences, \v{Z}itn\'a 25, CZ-11567 Praha 1, Czech
Republic ({\tt  krejci@math.cas.cz}).}
\and
Elisabetta Rocca
\thanks{Dipartimento di Matematica, Universit\`a di Milano,
Via Saldini 50, 20133 Milano, Italy ({\tt
elisabetta.rocca@unimi.it}). The work of ER was
supported by the FP7-IDEAS-ERC-StG Grant \#256872 (EntroPhase).}}

\date{}

\maketitle

\abstract{We propose an improved
model explaining the occurrence of high stresses due to the
difference in specific volumes during phase transitions between
water and ice. The unknowns of the resulting evolution problem are the absolute temperature,
the volume increment, and the liquid fraction. The main novelty here consists
in including the dependence of the specific heat and of the speed of sound
upon the phase. These additional nonlinearities bring
new mathematical difficulties which require new estimation techniques based on
Moser iteration. We establish the existence of a global solution to
the corresponding initial-boundary value problem, as well as lower and
upper bounds for the absolute temperature. Assuming constant heat conductivity,
we also prove uniqueness and continuous data dependence of the solution.}

\numberwithin{equation}{section}

\setcounter{section}{1}

\section*{Introduction}\label{intr}

In the present contribution we prove the well-posedness of an
initial-boundary value problem associated with the following
system coupling a quasi-linear parabolic internal energy balance
(for the absolute temperature $\teta$) with an
integro-differential equation for the relative volume increment
$U$,  and a differential inclusion ruling the evolution of the
phase variable $\chi$ as follows:
\begin{eqnarray}\nonumber
c(\chi)e_1(\theta)_t-
\dive(\kappa(\chi)\nabla\theta) &=& c'(\chi)\chi_t(f_1(\theta) -
e_1(\theta))\\ [2mm]\label{sys2}
&&\hspace{-28mm}
+\, \nu U_t^2 -\beta\theta U_t+
\gamma(\theta)\chi_t^2 - L\frac{\theta}{\theta_c}\chi_t\,,\\ [2mm]\label{sys1}
\nu U_t + \lambda(\chi) (U - \alpha (1-\chi)) - \beta (\theta - \theta_c)&=&
 \varrho_0 g (x_3 - \zeta_\Gamma)
- K_{\Gamma}(P_0(t)+U_\Omega(t))\,,\qquad\\ [2mm]\nonumber
&&\hspace{-40mm} -\gamma(\theta) \chi_t - \frac{\lambda'(\chi)}{2}
(U - \alpha(1-\chi))^2 - \alpha\lambda(\chi)
(U - \alpha(1-\chi))\\ [2mm]\label{sys3}
&&\hspace{-28mm}\in c'(\chi) \left(f_1(\theta) - f_1(\theta_c)\right)
+ L \left(1- \frac{\theta}{\theta_c}\right)+  \partial I(\chi)
\end{eqnarray}
with $U_\Omega(t) = \io U(x,t) \dd x$. In the previous paper
\cite{krsWil} we have already given a motivation and a complete
study of equilibria for this system, which  models the water
freezing in an elastic container, taking into account differences
in the specific volume, specific heat and speed of sound in the
solid and liquid phases. The derivation of the system from
physical principles and the meaning of the symbols will be
explained below in the next Section~\ref{equa}. Here, we describe
the mathematical difficulties and comment on previous results
related to this type of systems.

There is an abundant classical literature on phase transition
processes, see e.g. the monographs \cite{bs}, \cite{fremond},
\cite{visintin} and the references therein. It seems, however,
that only few publications take into account different mass
densities/specific volumes of the phases. In \cite{fr1}, the
authors proposed to interpret a phase transition process in terms
of a balance equation for macroscopic motions, and to include the
possibility of voids. Well-posedness of an initial-boundary value
problem associated with the resulting PDE system is proved there
and the case of two different densities $\varrho_1$ and
$\varrho_2$ for the two substances undergoing phase transitions
has been pursued in \cite{fr2}.

Here, we deal exclusively with physically measurable quantities.
All parameters have a clear physical meaning and the derivation is
carried out under the assumption that the displacements are small.
This enables us to state the system in Lagrangian coordinates (cf.
\cite{fr2} for a different approach to the subject).

The present model has been previously studied in \cite{krsbottle}
and \cite{krsgrav} under the assumption that the speed of sound
and the specific heat are the same in solid and in liquid. In
terms of the system \eqref{sys2}--\eqref{sys3}, this corresponds
to choose constant functions $\lambda(\chi)\equiv\lambda$ and
$c(\chi)\equiv c$. For this particular case, we have proved in
\cite{krsbottle} and \cite{krsgrav} the existence and uniqueness
of global solutions, as well as the convergence of the solutions
to equilibria. In reality, the specific heat in water is about the
double, while the speed of sound in water is less than one half of
the one in ice. The main goal of this contribution is to give a
well--posedness result for a boundary value problem associated
with (\ref{sys2})--(\ref{sys3}) including these dependences into
the model. The main result is stated in Section~\ref{mainres}. The
dependence of speed of sound and of the specific heat on the phase
is expressed in terms of additional nonlinearities in the
equations which have to be suitably handled. Moreover, here we
also generalize the results of \cite{krsbottle} and \cite{krsgrav}
allowing for non constant external pressure and temperature.
Finally, we proceed here with a different technique for the proof
of existence of solutions with respect to \cite{krsbottle} and
\cite{krsgrav}. Since the contraction argument does not work in
our situation, we discretize in time our problem (cf.
Subsection~\ref{approxi}), preparing thus necessary tools
for future numerical investigations on this model, and prove the
convergence of the scheme. The uniqueness and continuous
dependence of solution on the data is proved in
Section~\ref{proofuni} following the idea already exploited in
\cite{ckrsN} where we deal with a quasi-linear internal energy
balance equation coupled with a vectorial and nonlocal phase
dynamic. The main estimates are obtained here by means of the
energy inequality which still holds true at the discrete level
(cf. Subsection~\ref{esti}). Finally, it is worth noting that a
time dependent positive lower bound for the $\teta$-component of
the solution independent of the time step is established on the
time discrete approximation in Subsection~\ref{lowertetak}, while
we obtain a uniform in time upper bound on the solution $\teta$ by
means of a proper Moser estimate (cf. Subsection~\ref{concluexi}).

\section{Balance equations}\label{equa}

Referring to \cite{krsWil} for the complete deduction of the model,
we consider a liquid substance contained in a bounded connected
container $\Omega \subset \real^3$ with boundary of class $C^{1,1}$. The state variables are
the absolute temperature $\theta>0$, the displacement $\bfu \in \real^3$,
and the phase variable $\chi \in [0,1]$. The value  $\chi = 0$
means solid, $\chi = 1$ means liquid, $\chi \in (0,1)$ is a mixture
of the two.

We make the following modeling hypotheses.

\begin{itemize}
\item[{\bf (A1)}]
The displacements are small.
Therefore, we state the problem in {\em Lagrangian coordinates\/},
in which mass conservation is equivalent to the condition
of a constant mass density $\varrho_0>0$.
\item[{\bf (A2)}]
The substance is isotropic and compressible; the speed of sound and
the specific heat may depend on the phase $\chi$.
\item[{\bf (A3)}]
The evolution is slow, and we neglect shear viscosity and inertia effects.
\item[{\bf (A4)}]
We neglect shear stresses.
\item[{\bf (A5)}]
The liquid phase is the reference state, and the specific volume $V_i$ of the
solid phase is larger than the specific volume $V_w$ of the liquid phase.
\end{itemize}

We thus consider the evolution system
\begin{eqnarray}
-\dive \bfsi &=& \bff_{vol}\,,\label{bal1}\\[2mm]
\varrho_0 e_t + \dive \bfq &=& \bfsi\bfdd\bfe_t \,,\label{bal2} \\[2mm]
-\gamma_0(\theta)\chi_t &\in& \partial_\chi f\,,\label{bal3}
\end{eqnarray}
consisting of a the mechanical equilibrium equation \eqref{bal1}, energy conservation law \eqref{bal2},
and a phase dynamic equation \eqref{bal3}, where the coefficient $\gamma_0$ determines
the speed of the phase transition.
By {\bf (A4)}, the stress has the form $\bfsi = -p\,\bfde$ and the scalar quantity
\begin{equation}\label{press}
p := -\nu \bfe_t\bfdd\bfde - \lambda(\chi) (\bfe\bfdd\bfde - \alpha(1-\chi))
+ \beta (\theta - \theta_c)
\end{equation}
is the {\em pressure\/}. Here $\nu>0$ is a volume viscosity coefficient,
$\lambda(\chi)$ is the Lam\'e constant, which may depend on $\chi$ by virtue of
{\bf (A2)}, $\alpha=(V_i-V_w)/V_w$ is a positive phase expansion coefficient
by {\bf (A5)}, while $\beta$ is the
thermal expansion coefficient, which is assumed to be constant, and
$\bff_{vol}$ is a given volume force density (the gravity force)
\begin{equation}\label{grav}
\bff_{vol} = -\varrho_0 g\,\bfde_3\,,
\end{equation}
with standard gravity $g$ and vector $\bfde_3 = (0,0,1)$.

We denote by $e$ the specific internal energy,  $s$ is the specific entropy,
and $\bfq$ is the heat flux vector that we assume for simplicity in the form
\begin{equation}\label{flux}
\bfq = -\kappa(\chi)\nabla\theta
\end{equation}
with heat conductivity $\kappa(\chi) > 0$ depending possibly on $\chi$.

We assume the specific heat $c_V(\chi,\theta)$ in the form
\begin{equation}\label{speheat}
c_V(\chi,\theta) \ = \ c_0(\chi) c_1(\theta)\,.
\end{equation}
This is still a rough simplification, and further generalizations
are desirable. According to \cite[Chapter VI]{joo}
or \cite[Section 5]{made}, the purely caloric
parts $e_{cal}$ and $s_{cal}$ of the specific internal energy and specific entropy
are given by the formulas $e_{cal}(\chi,\theta) = c_0(\chi) e_1(\theta)$,
$s_{cal}(\chi,\theta) = c_0(\chi) s_1(\theta)$, with
\begin{equation}\label{speen}
e_1(\theta) = \int_0^\theta c_1(r)\dd r\,, \quad
s_1(\theta) = \int_0^\theta \frac{c_1(r)}{r}\dd r\,.
\end{equation}
Then, the specific free energy
$f = e - \theta s$ satisfies the conditions $\bfsi^e = \varrho_0 \partial_\bfe f$,
$s = - \partial_\theta f$.
With a prescribed constant latent heat $L_0$
and freezing point $\theta_c>0$ at standard atmospheric pressure $P_{stand}$,
the specific free energy $f$ necessarily has the form
\begin{eqnarray}\label{free}
f &=& c_0(\chi) f_1(\theta) + \frac{\lambda(\chi)}{2\varrho_0}
(\bfe\bfdd\bfde - \alpha(1-\chi))^2\\[2mm]\nonumber
&&
 -\, \frac{\beta}{\varrho_0}(\theta - \theta_c)\bfe\bfdd\bfde
 + L_0 \chi\left(1 - \frac{\theta}{\theta_c}\right)
 + \tilde f (\chi)\,,
\end{eqnarray}
where
\begin{equation}\label{f1}
f_1(\theta) \ = \ e_1(\theta) - \theta s_1(\theta) \ = \
\int_0^\theta c_1(r) \left(1 - \frac{\theta}{r}\right)\dd r\,,
\end{equation}
and $\tilde f$ is a arbitrary function of $\chi$
(integration ``constant'' with respect to $\theta$ and $\bfe$).
We choose $\tilde f$ so as to
ensure that the values of $\chi$ remain in the interval $[0,1]$, and that
the phase transition under standard pressure takes place
at temperature $\theta_c$. More specifically, we set
$$
\tilde f (\chi) \ = \  L_0 I(\chi) - c_0(\chi) f_1(\theta_c)\,.
$$
where $I$ is the indicator function of the interval $[0,1]$.

For specific entropy $s$ and
specific internal energy $e$ we obtain
\begin{eqnarray}\label{e2}
s &=& - \partial_\theta f = c_0(\chi) s_1(\theta)
+ \frac{\beta}{\varrho_0} \bfe\bfdd\bfde + \frac{L_0}{\theta_c}\chi\,,\\[2mm]\label{e3}%\nonumber
e &=& c_0(\chi)(e_1(\theta) - f_1(\theta_c))  +
\frac{\lambda(\chi)}{2\varrho_0} (\bfe\bfdd\bfde - \alpha(1-\chi))^2
 %\\[2mm]\label{e3} &&\qquad
+ \frac{\beta}{\varrho_0}\theta_c \bfe\bfdd\bfde + L_0(\chi + I(\chi)).\hspace*{10mm}
\end{eqnarray}

The equation for the phase $\chi$ is obtained by assuming that
$-\chi_t$ is proportional to $\partial_\chi f$
with proportionality coefficient (relaxation time) $\gamma_0(\theta)>0$,
where $\partial_\chi$ is the partial Clarke subdifferential with respect to $\chi$.

Then, the equilibrium equation (\ref{bal1}) can be rewritten in the form
$\nabla p = \bff_{vol}$, hence, as $\Omega$ is connected,
\begin{equation}\label{bala1}
p(x,t) = P(t) - \varrho_0 g\,x_3\,,
\end{equation}
where $P$ is a function of time only, which is to be determined. Recall
that in the reference state $\bfe\bfdd\bfde = \bfe_t\bfdd\bfde =0$,
$\chi=1$, and at standard pressure $P_{stand}$,
the freezing temperature is $\theta_c$. We thus see from (\ref{press})
that $P(t)$ is in fact the deviation from the standard pressure.
We assume also the external pressure
in the form $P_{ext} = P_{stand} + p_0$ with a given deviation $p_0(x,t)$.
The normal force acting on the boundary is
$(P(t) - \varrho_0 g\,x_3 -p_0) \bfn$,
where $\bfn$ denotes the unit outward normal vector.
We assume an elastic response of the boundary, and a heat transfer
proportional to the inner and outer temperature difference. On $\partial \Omega$,
we thus prescribe boundary conditions for $\bfu$ and $\theta$ in the form
\begin{eqnarray}\label{bcu}
(P(t) - \varrho_0 g\,x_3 -p_0(x, t))\bfn(x) &=& \bfk(x)\bfu(x,t) \,,\\[2mm]\label{bctheta}
\bfq(x,t)\cdot \bfn(x) &=& h(x) (\theta - \theta_\Gamma(x, t))
\end{eqnarray}
with a given symmetric positive definite matrix $\bfk(x)$ (elasticity of the boundary),
positive functions $h(x)$ (heat transfer coefficient), and
$\theta_\Gamma(x,t) >0$ (external temperature).
This enables us to find an explicit relation between $\dive \bfu$ and $P$.
Indeed, on $\partial \Omega$ we have by (\ref{bcu}) that
$\bfu\cdot\bfn = (P(t) - \varrho_0 g\,x_3 - p_0(x,t))\bfk^{-1}(x)\bfn(x)\cdot\bfn(x)$.
Assuming that $\bfk^{-1}\bfn\cdot\bfn$ belongs to $L^1(\partial \Omega)$, we set
\begin{equation}\label{bck}
\frac{1}{K_{\Gamma}} = \ipo \bfk^{-1}(x)\bfn(x)\cdot\bfn(x) \dd \sigma(x)\,,
\quad \zeta_\Gamma = K_\Gamma \ipo \bfk^{-1}(x)\bfn(x)\cdot\bfn(x)\, x_3 \dd \sigma(x)\,,
\end{equation}
and obtain by Gauss' Theorem that
\begin{equation}\label{bcu2}
U_\Omega(t) := \io \dive \bfu(x,t) \dd x = \frac{1}{K_{\Gamma}}(P(t)
- \varrho_0 g\, \zeta_\Gamma) - P_0(t)\,,
\end{equation}
where $P_0(t)= \int_{\partial\Omega} p_0(x,t)\bfk^{-1}(x)\bfn(x)\cdot\bfn(x) \dd\sigma(x)$.
Under the small strain hypothesis, the function $\dive \bfu$ describes the
local relative volume increment. Hence, Eq.~(\ref{bcu2}) establishes a linear
relation between the total relative volume increment $U_\Omega(t)$ and the relative pressure
$P(t) - p_0(x,t)$.
We have $\bfe\bfdd\bfde = \dive \bfu$, and thus the mechanical equilibrium equation
(\ref{bala1}), due to \eqref{press} and \eqref{bcu2}, reads
\begin{equation}\label{equau}
\nu \dive \bfu_t + \lambda(\chi) (\dive \bfu - \alpha(1-\chi)) - \beta
(\theta - \theta_c) + \varrho_0 g (\zeta_\Gamma - x_3) = -K_{\Gamma}(P_0(t)+U_\Omega(t))\,.
\end{equation}

As a consequence of (\ref{flux}), (\ref{free}), and (\ref{e3}),
the energy balance and
the phase relaxation equation in (\ref{bal2})--(\ref{bal3}) have the form
\begin{eqnarray} \nonumber
&&\hspace{-15mm} \varrho_0 c_0(\chi)e_1(\theta)_t - \dive(\kappa(\chi)\nabla\theta)
+ \varrho_0 c_0'(\chi)\chi_t(e_1(\theta) -
f_1(\theta))
\\[2mm]\label{equ1}
&=& \nu (\dive \bfu_t)^2 - \beta\theta \dive \bfu_t +
\varrho_0\gamma_0(\theta)\chi_t^2 - \varrho_0 L_0
\frac{\theta}{\theta_c}\chi_t\,, \\[2mm] \nonumber
&&\hspace{-15mm} -\varrho_0\gamma_0(\theta) \chi_t -
 \frac{\lambda'(\chi)}{2} (\dive\bfu - \alpha(1-\chi))^2 - \alpha\lambda(\chi)
(\dive\bfu - \alpha(1-\chi))
\\ [2mm]\label{equ4}
&\in& \varrho_0 c_1'(\chi) \left(f_1(\theta) -
f_1(\theta_c)\right) + \varrho_0 L_0
\left(\!1-\frac{\theta}{\theta_c}\!\right) +
 \partial I(\chi)\,.
\end{eqnarray}
For simplicity, we now set
\begin{equation}\label{zero}
U := \dive\bfu\,,\quad c(\chi) := \varrho_0 c_0(\chi)\,,
\quad \gamma(\theta) := \varrho_0\gamma_0(\theta)\,, \quad L := \varrho_0 L_0\,.
\end{equation}
Note that mathematically, the subdifferential $\partial I(\chi)$ is the same as
$\varrho_0 L_0 \partial I(\chi)$. The system thus reduces to the system (\ref{sys2})--(\ref{sys3}) of
three scalar equations -- one PDE and two ``ODEs'' for three
 unknown functions $\theta, \chi$, and $U$,
with boundary condition (\ref{bctheta}), (\ref{flux}). Assuming that a solution to \eqref{sys2}--\eqref{sys3}
is known with $U\in L^2(\Omega\times (0,T))$, we
find the vector function $\bfu$ by defining first $\Phi$ to be the solution to the Poisson
equation $\Delta\Phi = U$ with the Neumann boundary condition $\nabla\Phi\cdot \bfn =
(K_{\Gamma} U_\Omega(t) + \varrho_0 g(\zeta_\Gamma - x_3))
\bfk^{-1}(x)\bfn(x)\cdot\bfn(x)$.
With this $\Phi$, we find $\tilde \bfu$ as a solution to the problem
\begin{eqnarray}\label{sys4a}
\dive \tilde\bfu \ = \ 0\hspace{4.5mm} && \mbox{in }\ \Omega\times(0,T)\,,\\ \label{sys6a}
\left.
\begin{array}{rcl}
\tilde\bfu \cdot \bfn &=& 0 \\
(\tilde\bfu + \nabla\Phi - (K_{\Gamma} U_\Omega + \varrho_0 g(\zeta_\Gamma - x_3))
\bfk^{-1}\bfn)\times\bfn &=& 0
\end{array}
\right\}\hspace{2mm}
&& \mbox{on }\ \partial\Omega\times (0,T)\,,\qquad
\end{eqnarray}
and set $\bfu = \tilde\bfu + \nabla\Phi$. Then $\bfu$ satisfies a.e. in $\Omega$ the equation
$\dive \bfu = U$, together with the boundary condition
(\ref{bcu}), that is, $\bfu =
(K_{\Gamma} U_\Omega + \varrho_0 g(\zeta_\Gamma - x_3))
\bfk^{-1}\bfn$ on $\partial\Omega$.

For the solution to (\ref{sys4a})--(\ref{sys6a}), we refer to
\cite[Lemma 2.2]{gr} which states that
for each $\bfg\in H^{1/2}(\partial\Omega)^3$ satisfying $\ipo \bfg\cdot\bfn \dd \sigma(x) = 0$
there exists a function $\tilde\bfu\in H^1(\Omega)^3$, unique up to an additive
function $\bfv$ from the set $V$ of divergence-free $H^1(\Omega)$ functions vanishing on
$\partial\Omega$, such that $\dive \tilde\bfu = 0$ in $\Omega$, $\tilde\bfu = \bfg$
on $\partial\Omega$. In terms of the system (\ref{sys4a})--(\ref{sys6a}), it suffices to set
$\bfg =  ((\nabla\Phi -
(K_{\Gamma} U_\Omega + \varrho_0 g(\zeta_\Gamma - x_3))
 \bfk^{-1}\bfn)\times\bfn)\times\bfn$
and use the identity $(\bfb\times\bfn)\times\bfn = (\bfb\cdot\bfn)\,\bfn-\bfb$ for
every vector $\bfb$. Moreover, the estimate
\begin{equation}\label{turbu}
\inf_{\bfv \in V}\|\tilde \bfu + \bfv\|_{H^1(\Omega)} \le C\,\|\bfg\|_{H^{1/2}(\partial\Omega)}
\le \tilde C \|\Phi\|_{H^2(\Omega)}
\end{equation}
holds with some constants $C, \tilde C$. The required regularity is available here
by virtue of the assumption that $\Omega$ is of class $C^{1,1}$, provided $\bfk^{-1}$ belongs
to $H^{1/2}(\partial \Omega)$. Note that a weaker formulation
of problem (\ref{sys4a})--(\ref{sys6a}) can be found in \cite[Section 4]{ag}.

Due to our hypotheses {\bf (A3)}, {\bf (A4)}, we thus lose any control on
possible volume preserving turbulences $\bfv \in V$. This, however, has no influence
on the system (\ref{sys2})--(\ref{sys3}), which is the subject of our interest here.
Inequality (\ref{turbu}) shows that if $U$ is small in agreement
with hypothesis {\bf (A1)}, then also
$\bfv$ can be chosen in such a way that hypothesis {\bf (A1)},
interpreted in terms of $H^1$, is not violated.

%%%%%%%%%%%%%%%%%%%%%%%%%%%%%%%%%%%%%%%%%%%%%%%%%%%%%%%%%%%%%%%%%%%%%%%%%%%%%%%

\section{Energy and entropy}\label{ther}

In terms of the new variables $\theta, U, \chi$, the densities $\varrho_0 e, \varrho_0 s$ of
energy and entropy can be written~as
\begin{eqnarray}\label{energy}
\varrho_0 e &=& c(\chi)(e_1(\theta)-f_1(\theta_c))
+ \frac{\lambda(\chi)}{2}(U - \alpha(1-\chi))^2 % \\[2mm]\label{energy}&&
+ {\beta} \theta_c U + L(\chi + I(\chi))\,,\qquad \\[2mm]\label{entropy}
\varrho_0 s &=& c(\chi) s_1(\theta) + \frac{L}{\theta_c}\chi + {\beta} U \,.
\end{eqnarray}
The energy functional has to be supplemented with the boundary
energy term
\begin{equation}\label{enerb}
E_\Gamma(t) \ = \ \frac{K_{\Gamma}}{2}\left(U_{\Omega}(t) +P_0(t)+\frac{\varrho_0 g \zeta_\Gamma}{K_{\Gamma}}\right)^2,
\end{equation}
as well as with the gravity potential $-\varrho_0 g x_3 U$.
The energy and entropy balance equations now read
\begin{eqnarray}\nonumber
\frac{\dd}{\dd t} \left(\io \varrho_0 (e(x,t) - g x_3 U)\dd x
+ E_\Gamma(t) \right)
&=& \ipo h(x)(\theta_\Gamma(x,t) - \theta)\dd \sigma(x)\\[2mm]\label{princ1}
&&\hspace{-15mm}+\,K_\Gamma(P_0)_t(t) \left(U_{\Omega}(t) +P_0(t)
+\frac{\varrho_0 g \zeta_\Gamma}{K_{\Gamma}}\right),\\[2mm]\label{princ2}
\varrho_0 s_t + \dive \frac{\bfq}{\theta} &=& \frac{\kappa(\chi) |\nabla\theta|^2}{\theta^2}
+ \frac{\gamma(\theta)}{\theta} \chi_t^2 + \frac{\nu}{\theta}U_t^2 \ \ge\ 0,\qquad\\[2mm]\label{princ3}
\frac{\dd}{\dd t}\io \varrho_0 s(x,t)\dd x &=& \ipo \frac{h(x)}{\theta}(\theta_\Gamma(x,t) - \theta)
\dd \sigma(x) \\[1mm]\nonumber
&& \hspace{-15mm}+\, \io\left(\frac{\kappa(\chi) |\nabla\theta|^2}{\theta^2}
+ \frac{\gamma(\theta)}{\theta} \chi_t^2 + \frac{\nu}{\theta}U_t^2 \right)\dd x\,.
\end{eqnarray}
The entropy balance (\ref{princ2}) says that the entropy production
on the right hand side is nonnegative in agreement with the second principle
of thermodynamics. The system is not closed, and the energy supply
or the energy loss through the boundary is given by the right hand
side of (\ref{princ1}).

We prescribe the initial conditions
\begin{eqnarray}\label{ini2}
\theta(x,0) &=& \theta^0(x)\\ \label{ini1}
U(x,0) &=& U^0(x)\\ \label{ini3}
\chi(x,0) &=& \chi^0(x)
\end{eqnarray}
for $x \in \Omega$, and compute from (\ref{energy})--(\ref{entropy}) the corresponding
initial values $e^0$, $E_\Gamma^0$, and $s^0$ for specific energy,
boundary energy, and entropy, respectively.
Let $E^0$ and $S^0$ denote, respectively, $E^0 = \io \varrho_0 e^0(x)\dd x$,
$S^0=\io \vr_0s^0(x)\dd x$.
{}From the energy end entropy balance equations (\ref{princ1}), (\ref{princ3}),
we derive the following crucial (formal for the moment) balance equation for the
``extended'' energy $\varrho_0 (e - \bar \theta_\Gamma s)$, $\bar\theta_\Gamma$ being
a suitable positive constant:
\begin{eqnarray}\nonumber
&&\hspace{-10mm}\io\left(c(\chi)(e_1(\theta)- f_1(\theta_c)) +
\frac{\lambda(\chi)}{2}(U - \alpha(1-\chi))^2\right)(x,t)\dd x\\
\nonumber
&&+\, \io\left( \beta \theta_c U + L\chi - \varrho_0 g
x_3 U\right)(x,t)\dd x\\ \nonumber &&+\, \frac{K_{\Gamma}}{2}
\left(U_{\Omega}(t) + P_0(t)+\frac{\varrho_0
g\,\zeta_\Gamma}{K_{\Gamma}}\right)^2 \\ \nonumber
&&+\,\bar\theta_\Gamma\int_0^t\io \left(\frac{\kappa(\chi)
|\nabla\theta|^2}{\theta^2} + \frac{\gamma(\theta)}{\theta} \chi_t^2 +
\frac{\nu}{\theta}U_t^2  \right)(x,\xi)\dd x \dd\xi\\ \nonumber
&&+\, \int_0^t \ipo\frac{h(x)}{\theta}
(\theta - \theta_\Gamma(x,\xi))(\theta - \bar\theta_\Gamma)
\dd \sigma(x) \dd\xi\\[2mm]
\no
&=& E^0 +
E_\Gamma^0 - \bar\theta_\Gamma S^0+\bar \theta_\Gamma\io\left(c(\chi)
s_1(\theta) + \frac{L}{\theta_c}\chi + \beta U \right)(x,t)\dd
x\\[2mm]\label{crucial}
&&+\int_0^tK_\Gamma(P_0)_t(\xi)\left(U_{\Omega}(\xi) +P_0(\xi)+\frac{\varrho_0 g\zeta_\Gamma}{K_{\Gamma}}\right) \dd \xi\,.
\end{eqnarray}
We assume that both $c(\chi)$ and $\lambda(\chi)$ are bounded from above and
from below by positive constants. The growth of $s_1(\theta)$ is dominated
by $e_1(\theta)$ as a consequence of the inequality
$$
\frac{s_1(\theta) - s_1(\theta^*)}{e_1(\theta) - e_1(\theta^*)}
\le \frac{1}{\theta^*} \qquad \forall \theta> \theta^* > 0\,.
$$
We will use the relation \eqref{crucial} to get an upper bound for the solution
on the whole time interval $(0, \infty)$. From the identity
\begin{equation}\label{quadr}
\frac{1}{\theta} (\theta - a)(\theta - b) = \frac{1}{\theta} (\theta - \sqrt{ab})^2
- (\sqrt{b} - \sqrt{a})^2
\end{equation}
for all $\theta, a, b > 0$, it follows that we find
a constant $C>0$ independent of $t$ such that for all $t>0$ we have
\begin{eqnarray}\nonumber
&&\hspace{-16mm}\io\left(e_1(\theta) + U^2\right)(x,t)\dd x
+ \int_0^t\io \left(\frac{\kappa(\chi)|\nabla\theta|^2}{\theta^2}
+ \frac{\gamma(\theta)\chi_t^2}{\theta}  + \frac{\nu U_t^2}{\theta}  \right)(x,\xi)\dd x \dd\xi
\\ \label{esti1}
&&+\, \int_0^t \ipo\frac{h(x)}{\theta}
\left(\theta - \sqrt{\bar\theta_\Gamma\theta_\Gamma(x,\xi)}\right)^2\dd \sigma(x) \dd\xi \ \le \ C\,,
\end{eqnarray}
provided we assume that
\begin{equation}\label{hypoinf}
\int_0^\infty \ipo h(x)\left(\sqrt{\theta_\Gamma(x,t)} -
\sqrt{\bar\theta_\Gamma}\right)^2\dd \sigma(x) \dd t < \infty\,, \quad
\int_0^\infty|(P_0)_t(t)| \dd t < \infty\,.
\end{equation}

%%%%%%%%%%%%%%%%%%%%%%%%%%%%%%%%%%%%%%%%%%%%%%%%%%%%%%%%%%%%%%%%%

\section{Main results}
\label{mainres}

We construct the solution of (\ref{sys2})--(\ref{sys3}) by a combined truncation and
time discretization scheme. The method of proof is independent
of the actual values of the material constants, hence we choose for simplicity
\begin{equation}\label{const}
L= 2,\ \ \theta_c =\alpha = \beta= \nu =\vr_0= 1\,.
\end{equation}

We consider the following assumptions on the data.
\behy\label{hyp1} Assume that there exist positive constants $c_*$, $c^*$, $\underline{c}$,
$\bar c$, $\underline{\lambda}$, $\bar\lambda$, $\kappa_*$, $\lambda^*$, $\gamma_*$
such that
\begin{itemize}
\item[{\rm (i)}] $c$ convex, $c\in C^{1,1}([0,1])$, $0<c_*\leq c(z)$, $0<\underline{c}\leq
c'(z)\leq \bar c$, for all $z\in [0,1]$;
\item[{\rm (ii)}] $c_1\in
C^0(\real^+)$, $c_1(\theta) \ge c^*$ for $\theta\ge 1$,
$\lim_{\theta\to\infty} c_1(\theta)/\theta = \infty$, $e_1(\theta):=\int_0^\theta c_1(r) \dd r$,
$\int_0^1 c_1(r)/r\dd r<\infty$, $\int_0^1 c_1(r)/r^2 \dd r = \infty$;
\item[{\rm (iii)}] $\lambda$ convex,
$\lambda\in C^{1,1}([0,1])$, $0<\underline{\lambda}\leq \lambda(z)\leq\bar\lambda$,
$0 \ge \lambda'(z) \ge -\lambda^*$ for all $z\in [0,1]$;
\item[{\rm (iv)}] $\kappa\in C^{1,1}([0,1])$,
$0<\kappa_*\leq \kappa(z)$ for all $z\in [0,1]$;
\item[{\rm (v)}] $h\in L^\infty(\partial\Omega)$ is a non-negative function;
\item[{\rm (vi)}] $\gamma\in C^{0,1}(\real^+)$, $0<\gamma_*\leq\gamma(r)$ for all $r\in \real^+$.
\end{itemize}
\enhy

The liquid phase does not persist for very large temperatures and the behavior of
$c_1(\theta)$ as $\theta\to \infty$ thus cannot be experimentally verified.
We nevertheless believe that the growth condition (ii) in Hypothesis \ref{hyp1} is not completely
meaningless taking into account the fact that in the interval between $273$ and $373\,K$
(0--100${}^\circ$C), the function $c_1(\theta)$ is convex with a minimum at 35${}^\circ$C
\footnote{see
{\tt http://www.engineeringtoolbox.com/water-thermal-properties-d\_162.html.}}.

We introduce the following notation:
\begin{align}\label{symb1}
\AA(U,\chi,x,t) &:= \lambda(\chi)(U-1+\chi) + K_{\Gamma}(U_{\Omega}(t)+P_0(t))
+ g(\zeta_\Gamma-x_3)+ 1\,,\\ \label{symb2}
\BB(\chi, \theta) &:= c'(\chi)(f_1(\theta)-f_1(\theta_c))-2\theta\,,\\ \label{symb3}
\CC(U, \chi) &:= \frac{\lambda'(\chi)}{2}(U-1+\chi)^2+\lambda(\chi)(U-1+\chi) + 2\,.
\end{align}

System (\ref{sys2})--(\ref{sys3}) with boundary condition (\ref{bctheta}) then
can be written in the form
\begin{align}\nonumber
\io c(\chi)e_1(\theta)_t w(x)\dd x + \io \kappa(\chi)\nabla\theta \cdot \nabla w(x) \dd x &=
\,\ipo h(x)(\theta_\Gamma(x,t)-\theta) w(x) \dd\sigma(x)\\ \label{nequ1}
&\hspace{-62mm} - \io \big(U_t \AA(U,\chi,x,t)
+ \chi_t \big(\CC(U, \chi) + c'(\chi)(e_1(\theta)-f_1(\theta_c))\big)\big) w(x)\dd x\,,\\ \label{nequ2}
U_t  - \theta  &= - \AA(U,\chi,x,t)\,,\\ \label{nequ3}
\gamma(\theta)\chi_t + \BB(\chi, \theta)  +  \partial I(\chi) &\ni - \CC(U, \chi)\,,
\end{align}
where (\ref{nequ1}) is to be satisfied for all test functions $w \in W^{1,2}(\Omega)$
and a.e. $t>0$, while (\ref{nequ2})--(\ref{nequ3}) are supposed to hold a.e. in
the space-time cylinder that we denote $\Omega_T := \Omega\times (0,T)$ for $T>0$,
$\Omega_\infty := \Omega\times (0,\infty)$.

In this section we prove the following existence and uniqueness result.
\begin{theorem}\label{main}
Let Hypothesis \ref{hyp1} be satisfied, and let $\theta_\Gamma\in H^1(0,T; L^2(\partial\Omega))$
such that
$0< \theta_* \le \theta_\Gamma \le \theta^*$, and $P_0 \in W^{1,1}(0,T)$
be given functions. Let the initial conditions in \eqref{ini2}--\eqref{ini3} be such that
\begin{align}\nonumber
&\theta^0 \in W^{1,2}(\Omega)\cap L^\infty(\Omega)\,,
\quad 0<\theta_* \le \theta^0(x)\le \theta^*\quad \mbox{a.e.}\,,\\
\nonumber
&U^0, \chi^0 \in W^{1,2}(\Omega)\cap L^\infty(\Omega)\,,\quad 0 \le \chi^0(x) \le 1 \quad \mbox{a.e.}
\end{align}
Then there exists at least a solution $(\theta,U,\chi)$ to (\ref{nequ1})--(\ref{nequ3}),
(\ref{ini2})--(\ref{ini3}), and constants $\teta^\sharp(T) \ge \teta^\flat(T)>0$ such that
\begin{equation}\label{lowbou}
 \teta^\flat(T)\leq \teta(x, t)\leq \teta^\sharp(T)\quad \hbox{for a.e. }(x, t)\in \Omega_T\,,
\end{equation}
$\chi \in [0,1]$ a.e.,
$U, U_t,\chi_t \in L^\infty(\Omega_T)$,
$\theta_t\in L^2(\Omega_T)$,
$\nabla U, \nabla\chi, \nabla\theta \in L^\infty(0,T;L^2(\Omega))$.
If moreover condition \eqref{hypoinf} is satisfied, then the solution exists globally, and
 $\teta^\sharp(T)$ can be chosen independently of $T$.
Finally, if $\kappa(\theta)\equiv\bar\kappa\in \real^+$ is constant,
then the solution is unique, and its $L^2$-norm depends continuously on the data.
\end{theorem}

\begin{remark}\label{uni}
{\rm
Let us note that we could prove our existence result assuming that
$\kappa=\kappa(\theta,\chi)=k_1(\theta)k_2(\chi)$  with the same
techniques. Moreover, also uniqueness would hold true in case
$\kappa=\kappa_1(\theta)$ with an appropriate modification of the boundary condition
by means of the standard Kirchhoff transformation technique.
}
\end{remark}

%%%%%%%%%%%%%%%%%%%%%%%%%%%%%%%%%%%%%%%%%%%%%%%%%%%%%%%%%%%%%%%%%%%%%

\section{Existence proof}
\label{exiproof}

We proceed as follows: We truncate from above the functions depending on $\theta$ in
\eqref{nequ1}--\eqref{nequ3}, and discretize the system in time. For the discrete system,
we derive upper and lower bounds that enable us to let the time step tend to $0$ and prove
the existence of a solution to the truncated problem.
Finally, we prove
a time dependent lower bound and a uniform (in time and w.r.t. the truncation
parameters) upper bound on $\theta$, so that the truncation can be removed,
and this will conclude the proof of existence of solutions.

\subsection{Approximation and discrete energy estimate}
\label{approxi}

We introduce, for $\theta\in \real$, $R>0$, the functions
\begin{align}\label{defQR}
&Q_R(\theta)= \min\{\theta^+, B(R)\}, \quad B(R)=R^{1/2}(\min\{e_1(R), |f_1(R)|\})^{1/4},\\
\label{defic1ro}
&c_1^R(\theta)=c_1(Q_R(\theta)),\\
\label{defe1ro}
&e_1^R(\theta)=\int_0^\theta c_1^R(r)\, \dd r,\\
\label{defis1ro}
&s_1^R(\theta)=\int_0^\theta\frac{c_1^R(r)}{Q_R(r)}\dd r,\\
\label{deff1ro}
& f_1^R(\theta)=e_1^R(\theta)-Q_R(\theta)s_1^R(\theta)=\int_0^\theta c_1^R(r)\left(1-\frac{Q_R(\theta)}{Q_R(r)}\right)\dd r.
\end{align}

In the rest of the proof the following relations, which directly follow from the
above definitions and from Hypo.~\ref{hyp1}\,(ii), play an important role:
\begin{align}
\label{relaR1}
&\bullet\quad  \hbox{If } \theta \leq B(R)\hbox{ then } e_1^R(\theta)=e_1(\theta),
\ s_1^R(\theta)=s_1(\theta), \ f_1^R(\theta)=f_1(\theta),\\
\no
&\bullet \quad \hbox{If } \theta > B(R)\hbox{ then } e_1^R(\theta)=e_1(B(R))+c_1(B(R))(\theta-B(R)), \\
\label{relaR2}
&\qquad s_1^R(\theta)=s_1(B(R))+\frac{1}{B(R)}c_1(B(R))(\theta-B(R)),
\quad f_1^R(\theta)=f_1(B(R)),\\
\label{relaR3}
&\bullet\quad \hbox{If } \theta > R\hbox{ then } e_1^R(\theta)>e_1(R)>0,\quad f_1^R(\theta)<f_1(R)<0,\\
\label{relaR4}
&\bullet\quad \lim_{R\to\infty}\frac{e_1(R)}{R^2}=\lim_{R\to\infty}\frac{c_1(R)}{2R}=\infty,\\
\label{relaR5}
&\bullet\quad \lim_{R\to\infty}\frac{f_1(R)}{R^2}=-\lim_{R\to\infty}\frac{s_1(R)}{2R}=-\lim_{R\to\infty}
\frac{c_1(R)}{2R}=-\infty,\\
\label{relaR6}
&\bullet\quad \lim_{R\to\infty}\frac{|f_1(R)|}{B^2(R)}=\lim_{R\to\infty}\frac{e_1(R)}{B^2(R)}= \infty\,,\
\lim_{R\to\infty}\frac{B(R)}{R}=\infty.
\end{align}

We now introduce the time-discrete version of \eqref{nequ1}--\eqref{nequ3}. For an
arbitrary $n\in \nat$, we define the time step $\tau=T/n$. Choosing a constant
$c_R\in \real^+$ depending on $R$, which we specify below, we
look for a solution $\{(\theta_k, U_k, \chi_k)\}_{k=1}^n$ to the scheme
\begin{align}\no
&\hspace{-6mm} \frac{1}{\tau}\io c(\chik)\left(e_1^R(\tetak)-e_1^R(\tetakm)\right)w(x)\dd x
+\io\kappa(\chikm)\nabla\tetak \cdot \nabla w(x) \dd x
\\
\no
&\quad +\io c_R (\tetak\tetak^+-\tetakm\tetakm^+)w(x) \dd x +
\ipo h(x)(\tetak-\theta_{k\Gamma}) w(x) \dd \sigma(x)\\ \no
&= -\,
\io c'(\chik)\dchik(e_1^R(\tetakm)-f_1(\teta_c)) w(x)\dd x\\ \label{nequ1k}
&
\quad - \io \left(\dUk\AA_k(U_k,\chi_k,\chi_{k-1}, x)
+\dchik\CC_k(U_k, \chi_k, \chi_{k-1})\right)w(x)\dd x\,,\\
\label{nequ2k}
&\hspace{-6mm} \dUk-  Q_R(\tetakm) = - \AA_k(U_k,\chi_k,\chi_{k-1}, x)\,,\\
\label{nequ3k}
&\hspace{-6mm} \gamma(\tetakm)\dchik + \BB_k(\chi_k, \theta_{k-1}) + \partial I(\chik) \ni
- \CC_k(U_k, \chi_k, \chi_{k-1}),
\end{align}
where
\begin{align}\label{symbd1}
\AA_k(U_k,\chi_k,\chi_{k-1}, x) &:= \lambda(\chi_{k-1})(U_k-1+\chi_k) + K_{\Gamma}(U_{k\Omega}
+p_k)+ g(\zeta_\Gamma-x_3)+ 1\,,\\ \label{symbd2}
\BB_k(\chi_k, \theta_{k-1}) &:= c'(\chi_k)(f_1^R(\theta_{k-1})-f_1(\theta_c))-2Q_R(\tetakm)\,,\\ \label{symbd3}
\CC_k(U_k, \chi_k, \chi_{k-1}) &:= \frac{\lambda'(\chi_k)}{2}(U_k-1+\chi_k)^2
+\lambda(\chi_{k-1})(U_k-1+\chi_k) + 2\,,
\end{align}
with $\theta_{k\Gamma}=\theta_\Gamma(\cdot, k\tau)$, $p_k = P_0(k\tau)$,
$\Uko = U_\Omega(k\tau)$,
and with initial condition $(\theta_0, U_0, \chi_0)$ = $(\theta^0, U^0, \chi^0)$.
The symbol $\bdot{\vp}_k$ denotes the discrete time derivative $(\vp_{k}-\vp_{k-1})/\tau$
for a generic sequence $\{\vp_k\}$.
Eq.~(\ref{nequ1k}) has to be satisfied for all test functions $w \in W^{1,2}(\Omega)$,
while (\ref{nequ2k})--(\ref{nequ3k}) are supposed to hold a.e. in $\Omega$.

It is easy to see then that the latter semi-implicit scheme has a unique solution. Indeed,
at each time step, we assume that $\tetakm,\Ukm,\chikm$ are known, and find
$\Uk, \chik$ satisfying \eqref{nequ2k}--\eqref{nequ3k}.
For $\tau$ sufficiently small, \eqref{nequ2k}--\eqref{nequ3k} is an algebraic system for
$(\Uk, \chik)$ of the form $\Phi(\Uk,\chik) = Y_k$ with $\Phi:\real^2 \to \real^2$
strictly maximal monotone, hence it admits a unique solution. Finally, we insert
$\Uk$ and $\chik$ in \eqref{nequ1k} and solve the
resulting coercive elliptic equation, obtaining in that way the desired solution $(\tetak, \Uk, \chik)$.

Then, we note that the total energy balance still holds true for the discrete system.
Indeed, we take \eqref{nequ1k} with $w=1$, and denote
$E_k=e_1^R(\tetak)-f_1(\theta_c)$. We have
\begin{align}\no
&\hspace{-10mm} \frac{1}{\tau}\io c(\chik)(e_1^R(\tetak)-e_1^R(\tetakm))
+(\chik-\chikm)c'(\chik)(e_1^R(\tetakm)-f_1(\theta_c)) \dd x\\
\no
&=\frac{1}{\tau}\io c(\chik)(E_k-E_{k-1})+(\chik-\chikm)c'(\chik)E_{k-1}\dd x\\
\no
&=\frac{1}{\tau}\io c(\chik)E_k-(c(\chik)-(\chik-\chikm)c'(\chik))E_{k-1}\dd x\,,
\end{align}
and using the fact that $E_k\geq 0$ and that $c$ is convex (cf. Hypo.~\ref{hyp1} (i)), we get
\begin{align}\no
&\hspace{-10mm} \frac{1}{\tau}\io c(\chik)(e_1^R(\tetak)-e_1^R(\tetakm))
+(\chik-\chikm)c'(\chik)(e_1^R(\tetakm)-f_1(\theta_c)) \dd x\\
\label{B}
&\geq \frac{1}{\tau}\io \left(c(\chik)E_k-c(\chikm)E_{k-1}\right)\dd x\,.
\end{align}
Set now $S_k=\Uk+\chik-1$. Then, we obtain
\begin{align}\no
&\hspace{-10mm} \io (\dchik + \dUk)\lambda(\chikm)(\Uk+\chik -1)
+\dchik\frac{\lambda'(\chik)}{2}(\Uk+\chik-1)^2\dd x\\
\no
&=\frac{1}{\tau}\io S_k(S_k-S_{k-1})\lambda(\chikm)+\frac{1}{2}S_k^2(\chik-\chikm)\lambda'(\chik)\dd x\\
\no
&\geq\frac{1}{2\tau}\io(S_k^2-S_{k-1}^2)\lambda(\chikm)+S_k^2(\chik-\chikm)\lambda'(\chik)\dd x\\
\no
&=\frac{1}{2\tau}\io S_k^2\left(\lambda(\chikm)+\lambda'(\chik)(\chik-\chikm)\right)-\frac{1}{2\tau} S_{k-1}^2\lambda(\chikm)\dd x\,.
\end{align}
Using now the convexity of $\lambda$ (cf. Hypo.~\ref{hyp1} (iii)), we get
\begin{align}\no
&\hspace{-10mm} \io (\dchik+\dUk)\lambda(\chikm)(\Uk+\chi_k-1)
+\dchik\frac{\lambda'(\chik)}{2}(\Uk+\chik-1)^2\dd x\\
\label{C}
&\geq\frac{1}{2\tau}\io S_k^2\lambda(\chik)-S_{k-1}^2\lambda(\chikm)\dd x\,.
\end{align}
Hence, from \eqref{nequ1k}, using \eqref{B} and \eqref{C}, we obtain
\begin{align}\no
&\hspace{-10mm} \frac{1}{\tau}\io\left(c(\chik)E_k-c(\chikm)E_{k-1}
+\frac12 \left( S_k^2\lambda(\chik)-S_{k-1}^2\lambda(\chikm)\right) \right)\dd x\\
\no
&\quad +\bdot{U}_{k\Omega}\left(K_\Gamma(\Uko + p_k)+g\zeta_\Gamma\right)-\io g x_3 \dUk\dd x\\
\no
&\quad +\io\left(\dUk + 2\dchik
+c_R(\tetak\tetak^+-\tetakm\tetakm^+)\right)\dd x+\ipo h(x)(\tetak-\theta_{k\Gamma})\dd\sigma(x)\\
\label{D}
&\leq 0\,.
\end{align}
Summing now \eqref{D} over $k=1, \dots, m$, $1 \le m \le n$, we get
\begin{align}\no
&\hspace{-6mm} \io\left(c(\chi_{m})E_{m}+\frac12S_{m}^2\lambda(\chi_{m})
-g x_3U_{m}+U_{m} +\chi_{m}+c_R\tau \theta_{m}\theta_{m}^+\right)\dd x\\
\no
&\quad +\frac{K_\Gamma}{2}\left(U_{m\Omega} + p_m+\frac{g\zeta_\Gamma}{K_\Gamma}\right)^2
+\tau\sum_{k=1}^{m}\int_{\partial\Omega} h(x)(\tetak-\theta_{k\Gamma})\dd \sigma(x)\\
\no
&\leq\io\left(c(\chi_0)E_0+\frac12S_0^2\lambda(\chi_0)-g x_3U_0+U_0+\chi_0
+c_R \theta_0\theta_0^+\right)\dd x\\
\label{discren}
&\quad + \frac{K_\Gamma}{2}\left(U_{0\Omega} + p_0+\frac{g\zeta_\Gamma}{K_\Gamma}\right)^2
+ K_\Gamma \sum_{k=1}^{m} |p_k - p_{k-1}|
\max_{0\le k \le m} \left|U_{k\Omega} + p_k +\frac{g\zeta_\Gamma}{K_\Gamma}\right|\,.
\end{align}
Using the fact that $\tau\sum_{k=1}^n\ipo h(x)\theta_{k\Gamma}(x)\dd \sigma(x) \le C(T)$,
$\sum_{k=1}^n |p_k - p_{k-1}| \le C(T)$, with a constant $C(T)$ independent of $\tau$ and $R$,
we check that the left hand side of
\eqref{discren} is bounded independently of $\tau$ and $R$. Consequently, all terms in Eq.~\eqref{nequ2k} are
bounded by a multiple of $(1+B(R))$. Similarly, multiplying \eqref{nequ3k}
by $\dchik$ and using the fact that $\dchik \xi_k \ge 0$ for all $\xi_k \in \partial I(\chik)$,
we obtain the estimates
\begin{equation}\label{enerest}
\left.
\begin{array}{rcl}
|U_k|+|\dUk| &\leq& C(T)(1+B(R))\\
|\dchik| &\leq& C(T)(1+B(R)+B^2(R) + |f_1(B(R))|)
\end{array}
\right\}\quad \mbox{a.e.}
\end{equation}

\subsection{Lower bound for $\tetak$}
\label{lowertetak}

Here we derive a lower bound for the approximated absolute temperature $\tetak$.
We first rewrite \eqref{nequ1k} for $w\in W^{1, 2}(\Omega)$, $w\geq 0$ a.e.,
using \eqref{nequ2k}--\eqref{nequ3k}, in the form
\begin{align}\no
&\hspace{-6mm} \frac{1}{\tau}\io c(\chik)\left(e_1^R(\tetak)-e_1^R(\tetakm)\right)w(x)\dd x
+\io\kappa(\chikm)\nabla\tetak \cdot \nabla w(x) \dd x
\\ \no
&\quad +\io \left(c_R (\tetak\tetak^+-\tetakm\tetakm^+)
\right) w(x)\dd x + \ipo h(x)(\tetak - \theta_{k\Gamma}) w(x) \dd \sigma(x)
\\ \no
&\ge \io\dUk(\dUk - Q_R(\tetakm))w(x)\dd x\\ \label{nequlow}
&\quad + \io \dchik\big(\gamma(\tetakm)\dchik
+c'(\chik)(f_1^R(\tetakm)-e_1^R(\tetakm)) - 2Q_R(\tetakm)\big)w(x)\dd x\,,
\end{align}
where we have used again the fact that
$$
-\dchik \CC(\Uk,\chik,\chikm) \ \ge\ \dchik\big(\gamma(\tetakm)\dchik + \BB(\chik,\tetakm)\big)
$$
by definition of the subdifferential.
The right hand side of \eqref{nequlow} is bounded from below by
a negative multiple (depending on $R$) of $\tetakm\tetakm^+$.
We can now choose $c_R$ in \eqref{nequ1k} sufficiently large in order to get the following
inequality for all $w\in W^{1, 2}(\Omega)$, $w\geq 0$ a.e.:
\begin{align}\no
&\hspace{-6mm} \frac{1}{\tau}\io c(\chik)\left(e_1^R(\tetak)-e_1^R(\tetakm)\right)w(x)\dd x
+\io\kappa(\chikm)\nabla\tetak \cdot \nabla w(x) \dd x
\\
\label{low1}
&+\ipo h(x)(\tetak-\theta_{k\Gamma}) w(x) \dd \sigma(x)\geq -c_R\io \tetak\tetak^+\, w(x)\dd x\,.
\end{align}
We now compare this inequality with the constant decreasing sequence $\{v_k\}$ defined recurrently as
\begin{equation}\label{ode}
\frac{1}{\tau} c_*\left(e_1^R(v_k)-e_1^R(v_{k-1})\right)=-c_R v_k^2, \quad v_0:= \theta_*\,.
\end{equation}
We write \eqref{ode}, adding the zero term $-\dive(k(\chikm)\nabla v_k)$, in the form
\begin{align}\no
&\hspace{-6mm} \frac{1}{\tau}\io c_*\left(e_1^R(v_k)-e_1^R(v_{k-1})\right)w(x)\dd x
+\io\kappa(\chikm)\nabla v_k \cdot \nabla w(x) \dd x
\\
\label{low2}
& = -c_R\io v_k^2\, w(x)\dd x\,.
\end{align}
Subtracting \eqref{low1} from \eqref{low2} and
testing the difference by $w = H_\e(v_k-\tetak)$, where
 $H_\e$ is the regularization of the Heaviside function $H$,
\begin{equation}\label{heaviside}
H_\e(v)=\begin{cases}
0 &\quad\hbox{if } v\leq 0\\
v/\e & \quad \hbox{if } v\in (0,\e)\\
1 & \quad\hbox{if } v\geq \e
\end{cases}\,,
\end{equation}
we obtain, since $v_k < v_{k-1}$, that
\begin{align}\label{low3}
& \io  c(\chik)\left(\left(e_1^R(v_k)-e_1^R(v_{k-1})\right) - \left(e_1^R(\tetak)-e_1^R(\tetakm)\right)\right)
H_\e(v_k-\tetak)\dd x \le 0\,.
\end{align}
Assume that $\tetakm \ge v_{k-1}$ (this is true for $k=1$). For $\e \searrow 0$,
\eqref{low3} yields $\tetak\geq v_k$, and by induction we get
$\tetak\geq v_k> v_n$ for all $k=1,\dots, n$. By \eqref{ode}, we have
\[
e_1(v_k)-e_1(v_{k-1}) = - C\tau v_k^2
\]
with $C = c_R/c_*$.
Under Hypo.~\ref{hyp1} (ii), we have that the function $G(z)=-\int_z^{v_0}
\frac{c_1(s)}{s^2}\, ds =+\infty$.
Then, $G$ is increasing in $(0,v_0]$, $G(0+)=-\infty$, $G(v_0)=0$.
Moreover, by the Mean Value Theorem,
there exists $s_k\in [v_k, v_{k-1}]$ such that
\[
\frac{G(v_{k})-G(v_{k-1})}{e_1(v_k)-e_1(v_{k-1})}=\frac{G'(s_k)}{c_1(s_k)}=\frac{1}{s_k^2}\leq \frac{1}{v_k^2},
\]
hence $G(v_{k-1})-G(v_k)\leq C\tau$, that is, $G(v_n)\geq -Cn\tau$ and so
$\tetak \ge v_n \geq G^{-1}(-Cn\tau) = G^{-1}(-CT) =: \theta^\flat(T)$. This concludes the proof of
the lower bound for $\tetak$.

\subsection{Estimates}
\label{esti}

Now, we perform the estimates we need in order to pass to the limit as $\tau\searrow 0$
in (\ref{nequ1k})--(\ref{nequ3k}). The right hand side of \eqref{nequ1k} is bounded from above,
by virtue of \eqref{enerest}, by $C(T,R) (\tetakm + 1)$, where $C(T,R)$ is, here and in the sequel,
any sufficiently large constant depending only on $T$ and $R$, and independent of $k$ and $\tau$.
Testing \eqref{nequ1k} by $w = \tetak-\tetakm$, we obtain
\begin{align}\no
&\hspace{-6mm} \io\left(\frac{1}{\tau} c(\chik)\left(e_1^R(\tetak)
-e_1^R(\tetakm)\right)(\tetak-\tetakm)
+\kappa(\chikm)\nabla\tetak\nabla(\tetak-\tetakm)\right) \dd x\\
\no
&\quad +\io c_R(\tetak-\tetakm)^2(\tetak+\tetakm)\, \dd x
+ \int_{\partial\Omega}h(x)(\tetak-\teta_{k\Gamma})(\tetak-\tetakm)\, \dd \sigma(x)\\
\no
&\leq C(T,R)\io |\tetak-\tetakm|(\tetakm + 1) \dd x\,.
\end{align}
Using the lower bound for $\tetak$, and choosing $\chi_{-1}=\chi_0$, we get
\begin{align}\no
&\hspace{-3mm}\frac{1}{\tau C_1(T)}\io|\tetak-\tetakm|^2\dd x
+\io \left(\kappa(\chikm)|\nabla\tetak|^2-\kappa(\chi_{k-2})|\nabla\tetakm|^2\right)\dd x\\
\no
&
\quad+\ipo h(x)\left((\tetak-\teta_{k\Gamma})^2-(\tetakm-\teta_{(k-1)\Gamma})^2\right)
\dd \sigma(x)\\
\no
&\le \io \left(\kappa(\chikm)-\kappa(\chi_{k-2})\right)|\nabla\tetakm|^2\, \dd x
+\ipo h(x) |\teta_{k\Gamma}-\teta_{(k-1)\Gamma}|\,|\tetak-\teta_{k \Gamma}|\dd\sigma(x)\\
\no
&
\quad +\, \tau C(T,R)\io (\tetakm + 1)^2 \dd x\\
\no
& \le
\tau C(T,R)\left(\io \left(1 + |\tetakm|^2 + |\nabla\tetakm|^2\right)\dd x
+ \ipo h(x)\left(|\bdot\teta_{k\Gamma}|^2 + (\tetak-\teta_{k\Gamma})^2\right)\dd \sigma(x)\right),
\end{align}
where $C_1$ is a positive constant depending on $T$ but not on $\tau$.
The elementary inequality $\tetak^2 - \tetakm^2 \le \frac{1}{2\tau}(\tetak - \tetakm)^2
+ \frac{\tau}{2}(\tetak + \tetakm)^2$ enables us to rewrite
the above inequality in the form
\begin{equation}\label{q1}
q_k-q_{k-1}\leq\tau C(q_k+q_{k-1}+b_{k-1})\,,
\end{equation}
with
\begin{align}\no
q_k &= \io\left(\frac{2}{C_1(T)}\tetak^2 +\kappa(\chikm)|\nabla\tetak|^2\right)\dd x
+\ipo h(x)(\tetak-\teta_{k\Gamma})^2\dd \sigma(x)\,,\\ \no
b_k &= 1 + \ipo h(x)|\bdot\teta_{k\Gamma}|^2\dd \sigma(x)\,,
\end{align}
and $C := C(T,R)$. Inequality \eqref{q1} is equivalent to
\begin{equation}\label{q2}
q_k\leq\frac{1+\tau C}{1-\tau C} q_{k-1}+\frac{\tau C}{1-\tau C} b_{k-1}\,,
\end{equation}
which yields
\begin{equation}\label{q3}
q_k \ \leq\ \left(\frac{1+\tau C}{1-\tau C}\right)^k q_0
+\tau \frac{C}{1-\tau C} \sum_{j=0}^{k-1} b_j\left(\frac{1+\tau C}{1-\tau C}\right)^{k-1-j}
\ \leq\ {\rm e}^{3k\tau C}\left(q_0+\tau\sum_{j=0}^{k-1} b_j\right)
\end{equation}
holding true for $\tau\leq 1/(3C)$.
We conclude for all $m = 1, \dots, n$ that
\begin{equation}\label{esti2}
\frac{1}{\tau}\sum_{k=1}^{m}\io |\tetak-\tetakm|^2\, \dd x
+\io |\nabla\teta_{m}|^2\, \dd x
+\int_{\partial\Omega} h(x)(\teta_{m}-\teta_{m\Gamma})^2\, \dd \sigma(x)\leq C(T,R).
\end{equation}
Then, we introduce the piecewise constant and piecewise linear interpolants, for $t\in [(k-1)\tau, k\tau)$,
$k=1, \dots, n$, by the formula
\begin{equation}\label{interpol}
\underline\teta^{(\tau)}(x,t)=\tetakm(x), \ \ \bar\teta^{(\tau)}(x,t)=\tetak(x), \ \
\hat\teta^{(\tau)}(x,t)=\tetakm(x)+(t-(k-1)\tau)\,\bdot\tetak(x)\,,
\end{equation}
with a similar notation for $U$, $\chi$, $\theta_\Gamma$, and $P_0$. In particular, we set
\[
\hat e^{(\tau)}(x,t)= e_1^R(\tetakm(x))+ \frac{1}{\tau}(t-(k-1)\tau)
(e_1^R(\tetak(x)) - e_1^R(\tetakm(x)))\,.
\]
The previous estimates give immediately that
\[
\hat\teta_t^{(\tau)} \hbox{ bounded in } L^2(0,T; L^2(\Omega))\,,
\]
\[
\nabla \bar\teta^{(\tau)} \hbox{ bounded in } L^\infty(0,T; L^2(\Omega))\,,
\]
\[
\int_0^T\io\left(|\underline\teta^{(\tau)}-\hat\teta^{(\tau)}|^2 +
|\bar\teta^{(\tau)}-\hat\teta^{(\tau)}|^2\right)(x,t)\, \dd x \dd t\leq C(T,R)\tau^2\,.
\]
Letting $\tau$ tend to $0$ and passing to subsequences if necessary,
we get the convergences
\begin{equation}\label{conve}
\left.
\begin{array}{rcll}
\hat\teta^{(\tau)} &\to& \teta \quad &\hbox{strongly in } C^0([0,T]; L^2(\Omega))\,, \\
\underline\teta^{(\tau)} \ \to \ \teta\,, \ \bar\teta^{(\tau)} &\to& \teta \quad
&\hbox{strongly in } L^2(0,T; L^2(\Omega))\,, \\
\hat\teta_t^{(\tau)} &\to& \teta_t \quad
&\hbox{weakly in } L^2(0,T; L^2(\Omega))\,, \\
\nabla\bar\teta^{(\tau)} &\to& \nabla\teta \quad
&\hbox{weakly* in } L^\infty(0,T; L^2(\Omega))\,.
\end{array}
\right\}
\end{equation}
Now we estimate $\nabla\chik$ and  $\nabla\Uk$ as follows. From Eq.~\eqref{nequ2k} it follows that
\begin{align}\no
&(\dUk(x)-\dUk(y))(U_k(x)-U_k(y))\leq C(T,R)\Big( |U_k(x)-U_k(y)|^2+
|\chik(x)-\chik(y)|^2\\ \no
&\qquad +\,|\chikm(x)-\chikm(y)|^2+|\tetakm(x)-\tetakm(y)|^2+|x-y|^2\Big)
\end{align}
and, analogously, from \eqref{nequ3k}, we obtain
\begin{align}\no
& (\dchik(x)-\dchik(y))(\chi_k(x)-\chi_k(y))\leq C(T,R)\Big( |U_k(x)-U_k(y)|^2+
|\chik(x)-\chik(y)|^2\\ \no
&\qquad +\,|\chikm(x)-\chikm(y)|^2+|\tetakm(x)-\tetakm(y)|^2\Big)\,.
\end{align}
Summing up the two previous inequalities, we get
\begin{align}\no
&(U_k(x)-U_k(y))^2+(\chi_k(x)-\chi_k(y))^2\leq (U_{k-1}(x)-U_{k-1}(y))^2+(\chikm(x)-\chikm(y))^2\\ \no
&\qquad +\,\tau C(T,R)\Big((U_k(x)-U_k(y))^2+(\chi_k(x)-\chi_k(y))^2+(U_{k-1}(x)-U_{k-1}(y))^2\\ \no
&\qquad +\,(\chikm(x)-\chikm(y))^2+(\tetakm(x)-\tetakm(y))^2 +(x-y)^2\Big)\,.
\end{align}
We are again in the situation of Eq.~\eqref{q1}, with
$q_k=(U_k(x)-U_k(y))^2+(\chi_k(x)-\chi_k(y))^2$,
$b_k=(\tetakm(x)-\tetakm(y))^2+(x-y)^2$. Hence, by \eqref{q3}, we obtain
\begin{align} \no
&\hspace{-3mm}(U_k(x)-U_k(y))^2+(\chik(x)-\chik(y))^2\\ \no
&\leq C(T,R)\Big((U_0(x)-U_0(y))^2+(\chi_0(x)-\chi_0(y))^2
+\tau\sum_{j=0}^{k-1}(\teta_j(x)-\teta_j(y))^2+(x-y)^2\Big)\,.
\end{align}
Using now the previous estimate on $\nabla\bar\teta^{\tau}$, we get
\[
\nabla\hat U^{(\tau)}, \ \nabla\hat\chi^{(\tau)}\hbox{ bounded in } L^\infty(0,T; L^2(\Omega))\,.
\]
We already know that $\hat U_t^{(\tau)}, \hat\chi_t^{(\tau)}$ are bounded in $L^\infty(\Omega_T)$.
Furthermore,
\begin{eqnarray*}
\int_0^T\io|\bar U^{(\tau)}-\hat U^{(\tau)}|^2(x,t)\, \dd x \dd t &\leq& C(T,R) \tau^2\,,\\
\int_0^T\io\left(|\underline\chi^{(\tau)}-\hat \chi^{(\tau)}|^2 +
|\bar \chi^{(\tau)}-\hat \chi^{(\tau)}|^2\right)(x,t)\, \dd x \dd t &\leq& C(T,R) \tau^2,
\end{eqnarray*}
so that the convergences \eqref{conve} take place also for $U$ and $\chi$.
We now rewrite (\ref{nequ1k})--(\ref{nequ3k}) in terms of the functions
$\underline\teta^{(\tau)}, \bar\teta^{(\tau)}, \hat e^{(\tau)},
\underline\chi^{(\tau)}, \bar\chi^{(\tau)}, \hat\chi^{(\tau)}, \bar U^{(\tau)}$,
$\hat U^{(\tau)}, \bar\teta_\Gamma^{(\tau)}, \bar P_0^{(\tau)}$.
The above estimates allow us to pass to the limit as $\tau\searrow 0$ and obtain
a solution for the following truncated problem
\begin{align}\nonumber
\io c(\chi)e_1^R(\theta)_t w(x)\dd x + \io \kappa(\chi)\nabla\theta \cdot \nabla w(x) \dd x &=
\,\ipo h(x)(\theta_\Gamma(x,t)-\theta) w(x) \dd \sigma(x)\\ \label{nequ1R}
&\hspace{-62mm} - \io \big(U_t \AA(U,\chi,x,t)
+ \chi_t \big(\CC(U, \chi) + c'(\chi)(e_1^R(\theta)-f_1(\theta_c))\big)\big) w(x)\dd x\,,\\ \label{nequ2R}
U_t  - Q_R(\theta)  &= - \AA(U,\chi,x,t)\,,\\ \label{nequ3R}
\gamma(\theta)\chi_t + \BB_R(\chi, \theta)  +  \partial I(\chi) &\ni - \CC(U, \chi)\,,
\end{align}
where
\begin{align}\label{symb2R}
\BB_R(\chi, \theta) &:= c'(\chi)(f_1^R(\theta)-f_1(\theta_c))-2Q_R(\theta)\,,
\end{align}
$\AA$ and $\CC$ are defined in \eqref{symb1}, \eqref{symb3}, and
(\ref{nequ1R}) is to be satisfied for all test functions $w \in W^{1,2}(\Omega)$
and a.e. $t\in (0,T)$, while (\ref{nequ2R})--(\ref{nequ3R}) hold a.e. in
$\Omega_T$.

The next step consists in proving that $\theta$ remains uniformly bounded also from above
independently of $R$, so that the truncation does not become active if $R$ is sufficiently large.
The argument is based on the following counterpart of the extended energy balance \eqref{crucial},
\begin{eqnarray}\nonumber
&&\hspace{-10mm}\io\left(c(\chi)(e_1^R(\theta)- f_1(\theta_c)) +
\frac{\lambda(\chi)}{2}(U - 1+\chi)^2\right)(x,t)\dd x\\
\nonumber
&&+\, \io\left(U + 2\chi - g x_3 U\right)(x,t)\dd x +  \frac{K_{\Gamma}}{2}
\left(U_{\Omega}(t) + P_0(t)
+\frac{g \zeta_\Gamma}{K_{\Gamma}}\right)^2 \\ \nonumber
&&+\,\bar\theta_\Gamma\int_0^t\io \left(\frac{\kappa(\chi)|\nabla Q_R(\theta)|^2}{Q_R^2(\theta)}
+ \frac{\gamma(\theta) \chi_t^2}{Q_R(\theta)} +
\frac{U_t^2}{Q_R(\theta)} \right)(x,\xi)\dd x \dd\xi\\ \nonumber
&&+\, \int_0^t \ipo\frac{h(x)}{Q_R(\theta)}
(\theta - \theta_\Gamma(x,\xi))(Q_R(\theta) - \bar\theta_\Gamma)
\dd \sigma(x) \dd\xi\\[2mm]
\no
&=& E^0 +
E_\Gamma^0 - \bar\theta_\Gamma S^0+\bar \theta_\Gamma\io\left(c(\chi)
s_1^R(\theta) + 2\chi + U \right)(x,t)\dd x\\[2mm]\label{crucialR}
&&+\int_0^tK_\Gamma(P_0)_t(\xi)\left(U_{\Omega}(\xi) +P_0(\xi)
+\frac{g \zeta_\Gamma}{K_{\Gamma}}\right) \dd \xi\,,
\end{eqnarray}
which holds for every solution to \eqref{nequ1R}--\eqref{nequ3R} and every $t \in (0,T)$.

\subsection{Uniform upper bound for $\teta$}
\label{concluexi}

We choose $R$ large enough such that $B(R) > \theta^* \ge \bar\theta_\Gamma$. Then
$(\theta - \theta_\Gamma(x,\xi))(Q_R(\theta) - \bar\theta_\Gamma) \ge
(Q_R(\theta) - \theta_\Gamma(x,\xi))(Q_R(\theta) - \bar\theta_\Gamma)$. We may therefore
argue as at the end of Section \ref{ther} and obtain from \eqref{crucialR} for all $t \in (0,T)$ that
\begin{eqnarray}\nonumber
&&\hspace{-4mm}\io\left(e_1^R(\theta) + U^2\right)(x,t)\dd x
+ \int_0^t\io \left(\frac{\kappa(\chi)|\nabla Q_R(\theta)|^2 }{Q_R^2(\theta)}
+ \frac{\gamma(\theta)\chi_t^2}{Q_R(\theta)}  + \frac{U_t^2}{Q_R(\theta)}\right)(x,\xi)\dd x \dd\xi
\\ \no
&&\quad +\, \int_0^t \ipo\frac{h(x)}{Q_R(\theta)}
\left(Q_R(\theta) - \sqrt{\bar\theta_\Gamma\theta_\Gamma(x,\xi)}\right)^2\dd \sigma(x)\dd\xi\\ \label{esti3}
&& \le \, C \left(1 + \int_0^t \ipo h(x)\left(\sqrt{\theta_\Gamma(x,\xi)} -
\sqrt{\bar\theta_\Gamma}\right)^2\dd \sigma(x) \dd\xi + \int_0^t|(P_0)_t(\xi)| \dd\xi
\right).\qquad
\end{eqnarray}
In order to perform the Moser iteration scheme on $\teta$ as in \cite[Prop. 3.6]{krsgrav},
we need first to estimate $U$ and $U_t$ in terms of $\theta$. Rewriting \eqref{nequ2R} as
\[
U_t+\lambda(\chi) U= Q_R(\teta)+G(x,t)
\]
where, by virtue of \eqref{esti3}, $G(x,t)$ is bounded above by a positive constant
$G_0$. Denoting $\hat\lambda(x,t):=\int_0^t\lambda(\chi(x,s))\,\dd s$, we obtain
the formula
\[
U(x,t)={\rm e}^{-\hat\lambda(x,t)} U_0(x)+\int_0^t{\rm e}^{\hat\lambda(x,\xi)-\hat\lambda(x,t)}
\big(Q_R(\teta) + G\big)(x,\xi) \dd\xi\,.
\]
Using  Hypo.~\ref{hyp1} (iii), we get the estimate
\begin{equation}\label{esU}
|U(x,t)|\leq |U_0(x)|+\int_0^t{\rm e}^{\underline \lambda (t-\xi)}Q_R(\teta)(x, \xi)\dd \xi
+ \frac{G_0}{\underline \lambda}
\end{equation}
and
\begin{equation}\label{esUt}
|U_t(x,t)|\leq \bar\lambda|U(x,t)|+ Q_R(\teta)(x,t)+G_0\,.
\end{equation}

Now we are ready in order to start the Moser iteration scheme. Choose in \eqref{nequ1R} $w(x)=u^p$,
$u= \psi_R(\theta) := (Q_R(\theta) - R)^+,$
with any $p>1$ and with $R$ larger than the constants in Hypothesis \ref{hyp1}. Then
\[
\io c(\chi) (e_1^R(\teta))_t u^p\, \dd x +\frac{4p}{(p+1)^2}\io \kappa(\chi)|\nabla u^{\frac{p+1}{2}}|^2 \dd x
+\int_{\partial\Omega} h(x) u^{p+1} \dd s(x)
\]
\begin{equation}\label{moser1}
\le -\io \left(U_t \AA(U,\chi,x,t)-\chi_t\big(\CC(U, \chi)
- c'(\chi)(e_1^R(\theta)-f_1(\theta_c))\big)\right)u^p\dd x\,.
\end{equation}
Put $E_p^R(\teta)=\int_0^\teta c_1^R(r)\psi_R^p(r)\dd{\rm r}$.
Then, we can rewrite \eqref{moser1} as
\[
\io \left(c(\chi) E_p^R(\teta)\right)_t\dd x +\frac{4p}{(p+1)^2}\io \kappa(\chi)|\nabla u^{\frac{p+1}{2}}|^2 \dd x
+\int_{\partial\Omega} h(x) u^{p+1} \dd s(x)
\]
\[
\le - \io U_t  \AA(U,\chi,x,t) u^p\dd x
\]
\begin{equation}\label{moser2}
-\io\chi_t\left(\CC(U, \chi)u^p
+c'(\chi)\big((e_1^R(\theta)-f_1(\theta_c))u^p -E_p^R(\teta)\big)\right)\dd x\,.
\end{equation}
We now prove that the last integral in \eqref{moser2} is non-positive
if $R$ is sufficiently large.
First of all let us note that, if $\chi_t=0$ then
it vanishes. Hence, let us consider the case $\chi_t\neq 0$. Then, from \eqref{nequ3R}
it follows that
$\chi_t (\gamma(\theta)\chi_t + \BB_R(\chi, \theta) + \CC(U, \chi)) = 0$, hence
\[
\chi_t= -\frac{1}{\gamma(\teta)}\left(\BB_R(\chi, \theta) + \CC(U, \chi)\right).
\]
The last integral in \eqref{moser2} is of the form
$-\io \frac{1}{\gamma(\teta)}I_1\times I_2 \dd x$,
where
\begin{eqnarray*}
I_1 &:=& \chi_t = -\CC(U, \chi)
+2Q_R(\teta)+c'(\chi)(f_1(\teta_c)-f_1^R(\teta)),\\
I_2 &:=& \CC(U, \chi)u^p+c'(\chi)\big((e_1^R(\theta)-f_1(\theta_c))u^p-E_p^R(\teta)\big)
\\
&=& \big(\CC(U, \chi)-c'(\chi)f_1(\teta_c)\big)u^p + c'(\chi)\left(e_1^R(\teta)u^p-E_p^R(\teta)\right).
\end{eqnarray*}
We can now estimate from below the last term as follows
\[
e_1^R(\teta)u^p-E_p^R(\teta)=\int_0^{\teta}p e_1^R(r)u^{p-1}\dd {\rm r}\geq e_1(R) u^p.
\]
We have $I_2 = 0$ if $\theta\le R$, while for $\theta>R$ we have by Hypothesis \ref{hyp1}\,(i)
\begin{eqnarray*}
I_1 &\geq& -|\CC(U, \chi)|+2R + \underline c |f_1^R(R)-f_1(\teta_c)|,\\
I_2 &\geq& u^p\left(-|\CC(U, \chi)|+ \underline c f_1(\teta_c) + \underline c e_1^R(R)\right).
\end{eqnarray*}
By virtue of \eqref{enerest}, we have
$|\CC(U, \chi)|\leq C(B^2(R)+1)$. Referring to \eqref{relaR6},
we conclude that there exists $R_0>1$ larger than all constants in Hypothesis \ref{hyp1}
such that for $R\geq R_0$ we have in \eqref{moser2}
\[
-\io\chi_t\left(\CC(U, \chi)u^p
+c'(\chi)\big((e_1^R(\theta)-f_1(\theta_c))u^p -E_p^R(\teta)\big)\right)\dd x \le 0\,.
\]

Let us fix now $R>R_0$ and continue the Moser estimate, rewriting \eqref{moser2} as follows
\[
\io \left(c(\chi) E_p^R(\teta)\right)_t\, \dd x +\frac{4p}{(p+1)^2}\io \kappa(\chi)|\nabla u^{\frac{p+1}{2}}|^2 \dd x
+\int_{\partial\Omega} h(x) u^{p+1} \dd s(x)
\]
\begin{equation}\label{moser3}
\leq
- \io U_t  \AA(U,\chi,x,t) u^p\dd x\,.
\end{equation}
We have $c(\chi)E_P^R(u)\geq \frac{c_*c^*}{p+1} u^{p+1}$, $\kappa(\chi) \ge \kappa_*$.
Integrating \eqref{moser3} in time, we obtain, using Hypo.~\ref{hyp1}\,(i),(ii),(iv), as well as
the estimates \eqref{esU}--\eqref{esUt} and the fact that $u(x,0) \equiv 0$, that
\[
\frac{c_*c^*}{p+1} \io u^{p+1}(x,t) \dd x +\frac{4p \kappa_*}{(p+1)^2}
\int_0^t\io \left|\nabla u^{\frac{p+1}{2}}\right|^2(x, \xi) \dd x\dd \xi
\]
\[
\leq\int_0^t \io u^p(x,\xi)\left(1+ Q_R(\teta(x, \xi))
+\int_0^\xi e^{-\underline\lambda (\xi-\eta)}Q_R(\teta(x,\eta))\dd\eta\right)r(x,\xi)\dd x\dd \xi\,.
\]
The function $r(x,t) = C \AA(U,\chi,x,t)$, where $C$ is a suitable constant, has norm in
$L^\infty (0,T; L^2(\Omega))$  bounded independently of $R$ by virtue of \eqref{esti3}.
Note that $Q_R(\theta) \le u + R$. Hence, the function $v:= u/R$ satisfies for all $p>1$
the inequality
\[
\frac{c_*c^*}{p+1} \io v^{p+1}(x,t) \dd x +\frac{4p \kappa_*}{(p+1)^2}
\int_0^t\io \left|\nabla v^{\frac{p+1}{2}}\right|^2(x, \xi) \dd x\dd \xi
\]
\[
\leq\int_0^t \io v^p(x,\xi)\left(1+ v(x, \xi)
+\int_0^\xi e^{-\underline\lambda (\xi-\eta)}v(x,\eta)\dd\eta\right)r(x,\xi)\dd x\dd \xi.
\]
The argument of \cite[Prop. 4.5]{krsbottle} yields $\|v\|_{L^\infty(\Omega_T)}\leq \bar C$
with a constant $\bar C$ independent of $R$ and $T$. Consequently,
\[
\|Q_R(\theta)\|_{L^\infty(\Omega_T)}\leq (1+ \bar C) R.
\]
Choosing $R$ sufficiently large such that $B(R) > (1+ \bar C) R$, we can
remove the truncation from (\ref{nequ1R})--(\ref{nequ3R}), concluding
in this way the proof of existence of a bounded solution to (\ref{nequ1})--(\ref{nequ3}).
If moreover \eqref{hypoinf} holds, then $r \in L^\infty(0,\infty; L^2(\Omega))$, and the
upper bound holds globally in $\Omega_\infty$. Indeed,
 the lower bound for $\teta$ in Subsection \ref{lowertetak}
is independent of the time step $\tau$ and is preserved when $\tau \searrow 0$.

%%%%%%%%%%%%%%%%%%%%%%%%%%%%%%%%%%%%%%%%%%%%%%%%%%%%%%%%%%%%%%%%%%%%%%%%%

\section{Uniqueness and continuous data dependence}
\label{proofuni}

In this Section, we prove uniqueness and continuous data dependence
of solutions under the more restrictive assumption that
$\kappa(r)=\bar\kappa\in \real^+$ for all $r\in \real^+$.

In what follows,
we denote by $R_0, R_1, R_2, \dots$ suitable constants that possibly depend on  $T$,
but not on the solutions. We first rewrite Eq.~\eqref{nequ1} in the form
\begin{align}\nonumber
\io(c(\chi)(e_1(\theta)-f_1(\theta_c)))_t w(x)\dd x+\io \bar\kappa\nabla\theta \cdot \nabla w(x)\dd x
&+ \,\ipo h(x)(\theta-\theta_\Gamma) w(x) \dd \sigma(x)\\ \label{uniequ1}
&\hspace{-52mm} = - \io \big(U_t \AA(U,\chi,x,t) + \chi_t \CC(U, \chi)\big) w(x)\dd x\,,
\end{align}
and denote
$\hat\theta = \theta_1- \theta_2$,
$\hat\chi = \chi_1-\chi_2$, $\hat\chi_0 = \chi_{01}-\chi_{02}$,
$\hat\theta_0 = \theta_{01}- \theta_{02}$, $\hat\theta_{\Gamma} = \theta_{\Gamma 1}
- \theta_{\Gamma 2}$, $\hat U=U_1-U_2$, $\hat U_0=U_{01}-U_{02}$,
$\hat\Theta(x,t) = \int_0^t \hat\theta(x,\tau) \dd\tau$, $\hat\Theta_\Gamma(x,t) =
\int_0^t \hat\theta_\Gamma(x,\tau) \dd\tau$, $\hat P_0 = P_{01}-P_{02}$.
Within the range $\teta^\flat(T) \le  \teta \le \teta^\sharp(T)$
and $\chi \in [0,1]$, $|\chi_t| \le C$ of admissible values for the solutions, and,
thanks to Hypo.\,\ref{hyp1}, all nonlinearities
in (\ref{nequ1})--(\ref{nequ3}) are Lipschitz continuous.
We integrate the difference of the two equations (\ref{uniequ1}),
written for $(\theta_1, U_1, \chi_1)$ and $(\theta_2, U_2, \chi_2)$, from $0$ to $t$, and test
by $w = \theta_1-\teta_2$. This yields
\begin{align}\nonumber
&\io|\hat\theta(x,t)|^2\,\dd x +
\frac{\dd}{\dd t} \left(R_0\io |\nabla\hat\Theta(x,t)|^2 \dd x
+ R_1\ipo h(x) (\hat\Theta - \hat\Theta_\Gamma)^2(x,t) \dd \sigma(x)\right)\\
\nonumber
&\ \leq R_2\Big(\|\hat\theta_0\|_{L^2(\Omega)}^2 + \itt|\hat P_0(\xi)|^2\dd\xi
+ \ipo h(x)|\hat\Theta - \hat\Theta_\Gamma|\,
|\hat\theta_\Gamma|(x,t) \dd \sigma(x)\\
\label{euni1}
&
+\io \Big(\itt \big(|\hat\chi_t(x,\xi)|+|\hat U_t(x,\xi)|+|\hat \chi(x,\xi)|+ |\hat U(x,\xi)|\big)\dd\xi\Big)^2\dd x \Big).
\end{align}
Repeating the argument of \cite[Proposition 4.3]{krsbottle} or \cite[Proposition 3.4]{ckrsN}
about the $L^1$-Lipschitz continuity of solution operators to gradient flows,
we obtain for the solutions to (\ref{nequ2})--(\ref{nequ3}) for a.e. $(x,t) \in \Omega_T$
the estimate
\begin{align}\label{uni1}
&\itt (|\hat\chi_t(x,\tau)| + |\hat U_t(x,\tau)|)(x,\xi)\dd\xi+|\hat\chi(x,t)|
+|\hat U(x,t)|\\
\no
&
\leq R_3 \left(\itt|\hat P_0(\xi)|\dd\xi+ |\hat\chi_0(x)|+|\hat U_0(x)|
+\itt\left(|\hat\theta(x,\xi)| + \io |\hat\theta(y,\xi)|\dd y \right)\dd\xi\right).
\end{align}
Integrating \eqref{euni1} from $0$ to $t$ and using \eqref{uni1} together with
Gronwall's argument, we obtain for each $t\in [0,T]$ the estimate
\begin{align}\no
&\hspace{-10mm}\int_0^t \io |\hat\theta(x,\xi)|^2 \dd x\dd\xi + \io|\hat\chi(x,t)|^2\dd x + \io|\hat U(x,t)|^2\dd x\\
\no
&\le R_4 \Big(\|\hat\theta_0\|_{L^2(\Omega)}^2 + \|\hat\chi_0\|_{L^2(\Omega)}^2 +\|\hat U_0\|_{L^2(\Omega)}^2 +\itt|\hat P_0(\xi)|^2\dd\xi\\
\label{uni8}
&
\hspace{10mm}+ \int_0^t\int_{\partial\Omega} h(x)\hat\theta_\Gamma^2(x,\xi) \dd \sigma(x)\dd\xi
\Big).
\end{align}
This concludes the proof of uniqueness
 of solutions and of Theorem~\ref{main}.

%\printindex
\end{document}